\documentclass[12pt]{amsart}
\usepackage[top=1in, bottom=1in, left=1in, right=1in]{geometry}
\usepackage{dsfont}
\usepackage{amsmath}
\usepackage{amsthm}
\usepackage{amssymb}
\usepackage{mathtools}
\usepackage{times}
\usepackage[utf8]{inputenc}
\usepackage{indentfirst}  
\usepackage{graphicx}
\usepackage{color}
\usepackage[colorlinks=true]{hyperref}
\hypersetup{colorlinks=true, citecolor=green, linkcolor=green, filecolor=magenta, urlcolor=cyan}
\usepackage[shortlabels]{enumitem}
\usepackage{mdframed}
\usepackage{tcolorbox}
\newtcolorbox{mybbox}{colback=blue!5!white,
colframe=blue!75!black}
\newtcolorbox{mygbox}{colback=green!5!white,
colframe=green!75!black}
\newtcolorbox{myrbox}{colback=red!5!white,
colframe=red!75!black}

\numberwithin{equation}{section}

\newtheorem{thm}{Theorem}[section]
\newtheorem{prop}[thm]{Proposition}
\newtheorem{lem}[thm]{Lemma}
\newtheorem{cor}[thm]{Corollary}
\theoremstyle{definition}

\theoremstyle{definition}

\newenvironment{feqn*}{\begin{mdframed}\begin{equation*}}{\vspace{1mm}
\end{equation*}\end{mdframed}}

\theoremstyle{remark}
\newtheorem{remark}{Remark}[section]
\newtheorem*{remark*}{Remark}

\newcommand{\N}{\mathbb{N}}
\newcommand{\Z}{\mathbb{Z}}

\newcommand{\R}{\mathbb{R}}
\newcommand{\C}{\mathbb{C}}

\newcommand{\CA}{\mathcal{A}}

\newcommand{\CC}{\mathcal{C}}

\newcommand{\CP}{\mathcal{P}}

\newcommand{\CR}{\mathcal{R}}

\newcommand{\m}{\mathfrak{m}}

\newcommand{\dee}{\mathrm{d}}

\newcommand{\abs}[1]{\lvert #1 \rvert}

\newcommand{\flbgg}[1]{\bigg\lfloor #1 \bigg\rfloor}

\newcommand{\bs}\boldsymbol{}
\renewcommand{\geq}{\geqslant}
\renewcommand{\leq}{\leqslant}
\renewcommand{\hat}{\widehat}

\newcommand{\eps}{\varepsilon}
\renewcommand{\Re}{{\rm Re}}

\renewcommand{\mod}[1]{\,({\rm mod}\,#1)}

\definecolor{blue}{rgb}{.2,.6,.75}
\definecolor{green}{rgb}{.4,.7,.4}
\definecolor{red}{rgb}{1,0,0}

\begin{document}

\title[A Pretentious Proof of Linnik's Estimate for Primes in Arith. Prog.]{A Pretentious Proof of Linnik's Estimate for Primes in Arithmetic Progressions}

\author{Stelios Sachpazis}

\address{D\'epartement de math\'ematiques et de statistique\\
Universit\'e de Montr\'eal\\
CP 6128 succ. Centre-Ville\\
Montr\'eal, QC H3C 3J7\\
Canada}

\email{stelios.sachpazis@umontreal.ca}


\date{\today}

\begin{abstract}
In the present paper, we adopt a pretentious approach and prove a strongly uniform estimate for the sums of the von Mangoldt function $\Lambda$ on arithmetic progressions. This estimate is analogous to an estimate that Linnik established in his attempt to prove his celebrated theorem concerning the size of the smallest prime number of an arithmetic progression. Our work builds on ideas coming from the pretentious large sieve of Granville, Harper and Soundararajan and it also borrows insights from the treatment of Koukoulopoulos on multiplicative functions with small averages.
\end{abstract}

\maketitle

\section{Introduction}

Let $q$ be a positive integer and $a\in(\Z/q\Z)^*$. We define $p(q,a)$ as the least prime occurring in the arithmetic progression $\{qn+a, n\in\N\}$. Under the Riemann Hypothesis for the $L$-Dirichlet series of characters $\chi$ mod $q$, it can be proved that 

\vspace{-2mm}

$$p(q,a)\ll(\phi(q)\log q)^2,$$

\vspace{2mm}

\noindent
where $\phi$ denotes the Euler totient function. Unconditionally, using the Siegel-Walfisz Theorem, one can show that

\vspace{-2mm}

$$p(q,a)\ll_{\eps}\exp(q^{\eps})$$

\vspace{2mm}

\noindent
for any $\eps>0$. Moreover, if there are no Siegel zeros, this bound reduces to

\vspace{-1mm}

$$p(q,a)\ll q^{c\log q}$$

\vspace{2mm}

\noindent
for some $c>0$. So, in 1944, it came as a big surprise when Linnik \cite{lin1, lin2} proved that there exist positive universal constants $C$ and $L$ such that

\vspace{-2mm}

$$p(q,a)\leq Cq^L$$

\vspace{2mm}

\noindent
for any choice of $q\in\N$ and $a\in(\Z/q\Z)^*$. One of the steps that led Linnik to this breakthrough, now known as Linnik's theorem, was an analogue of the following estimate. We will be referring to it as Linnik’s estimate since it is closely related to Linnik’s theorem.

Assume that $x\geq q^2$. If there exists an exceptional character $\chi$ corresponding to a Siegel zero, then

\vspace{-2mm}

\begin{equation}\label{linnest}
\sum_{\substack{n\leq x\\n\equiv a\mod{q}}}\Lambda(n)=\frac{x}{\phi(q)}+\frac{\chi(a)}{\phi(q)}\sum_{n\leq x}\Lambda(n)\chi(n)+O\bigg(\frac{x^{1-c_1/\log(2q)}}{\phi(q)}+\frac{xe^{-c_2\sqrt{\log x}}}{\phi(q)}\bigg),
\end{equation}

\noindent
where $c_1$ and $c_2$ are two positive absolute constants. If the exceptional character $\chi$ does not exist, then (\ref{linnest}) remains true, but the character sum of the right-hand side is omitted. 

Linnik's work on $p(q,a)$ was later simplified, as was done by Bombieri in \cite{bomb}, but the new proofs, including Linnik's original proof, relied in one form or another on the following three ingredients.

\vspace{2mm}

\begin{itemize}
	\item The classical zero-free region;
	\item A log-free zero-density estimate;
	\item The exceptional zero repulsion, also known as the Deuring-Heilbronn phenomenon, stating that it is possible to enlarge the classical zero-free region when it contains a Siegel zero.
\end{itemize}

\vspace{2mm}

\noindent
Proofs that make use of these three principles can be found in modern treatments, like the one which is presented in \cite[Chapter 18]{iwko}.

Nonetheless, in the last years, some proofs that avoid these tools have appeared. For example, in 2002, such a proof was developed by Elliot in \cite{el}. Later, in 2016, Granville, Harper and Soundararajan \cite{ghal} studied pretentiously the distribution of multiplicative functions on arithmetic progressions and as a consequence of their general results they were able to show a weaker form of (\ref{linnest}). In turn, this served as the stepping stone for another new proof of Linnik's theorem which circumvented the combination of the three aforementioned principles. Another pretentious proof of Linnik's theorem is presented in \cite[Chapter 27]{dimb} and a basic element of the proof is a flexible variant of (\ref{linnest}) where every prime is weighted with $1/p$ instead of $\log p$. Even though the alternative approaches recover Linnik's theorem, they do not provide a pretentious proof for an estimate of the same quantitative strength as (\ref{linnest}).

To this end, in this paper, we apply pretentious methods and prove Linnik's estimate (\ref{linnest}) with a refined error term where a Korobov-Vinogradov-type term replaces $\exp\{-c_2\sqrt{\log x}\}$. In particular, we prove the following theorem.

\begin{thm}\label{main}
Let $q\geq 1$ be an integer and consider a real number $x\geq q^2$. For any character $\chi$ mod $q$, we set 

$$L_q(1,\chi)=L(1,\chi)\prod_{p\leq q}(1-\chi(p)/p).$$

\vspace{1mm}

\noindent
We also define $\CR_q$ as the set of real, non-principal characters mod $q$ and we take a character $\psi$ such that $L_q(1,\psi)=\min_{\chi\in\CR_q}L_q(1,\chi)$. Then, for any $a\in(\Z/q\Z)^*$, we have that

\begin{align*}
\sum_{\substack{n\leq x\\n\equiv a\mod{q}}}\Lambda(n)=\frac{x}{\phi(q)}+\frac{\psi(a)}{\phi(q)}\sum_{n\leq x}\psi(n)\Lambda(n)+O\bigg(\frac{x^{1-C_1/\log(2q)}}{\phi(q)}+\frac{xe^{-C_2(\log x)^{3/5}(\log \log x)^{-3/5}}}{\phi(q)}\bigg),
\end{align*}

\vspace{3mm}

\noindent
where $C_1$ and $C_2$ are two positive absolute constants.
\end{thm}

We can bound the sum $\sum_{n\leq x}\psi(n)\Lambda(n)$ in Theorem \ref{main} by referring to Theorem 1.6(a) of \cite{oldk}. According to this theorem, there exist positive constants $c'$ and $c''$ such that

\begin{equation}\label{forrem}
\sum_{n\leq x}\psi(n)\Lambda(n)\ll x^{1-c'L_q(1,\psi)/\log(2q)}+xe^{-c''\sqrt{\log x}},
\end{equation}

\noindent
for all $x\geq q^2$. The same theorem also provides information about the size of the quantity $L_q(1,\psi)$ that is involved in the bound above. It is known that there exists a constant $\delta\in(0,1)$ such that $L(\cdot,\psi)$ has at most one zero $\beta$ in $[1-\delta/\log q,1)$. If such a zero does not exist, we set $\beta=1-\delta/\log q$. In either case, Theorem 1.6(a) of \cite{oldk} claims that $L_q(1,\psi)\asymp (1-\beta)\log q$. Therefore, we arrive at the following consequence of Theorem \ref{main}.

\begin{cor}\label{1stc}
Under the considerations and assumptions of Theorem \ref{main}, we have that

$$\sum_{\substack{n\leq x\\n\equiv a\mod{q}}}\Lambda(n)=\frac{x}{\phi(q)}+O\bigg(\frac{x^{1-\alpha_1L_q(1,\psi)/\log(2q)}}{\phi(q)}+\frac{xe^{-\alpha_2\sqrt{\log x}}}{\phi(q)}\bigg),$$

\vspace{1mm}

\noindent
for some positive constants $\alpha_1$ and $\alpha_2$. Moreover, there exists a constant $\delta\in(0,1)$ such that $L(\cdot,\psi)$ has at most one zero $\beta$ in $[1-\delta/\log q,1)$. If such a zero does not exist, we put $\beta=1-\delta/\log q$ and in any case, $L_q(1,\psi)\asymp (1-\beta)\log q$.
\end{cor}

\begin{remark}
The term $x^{1-c'L_q(1,\psi)/\log(2q)}$ in estimate (\ref{forrem}) leads to the first fraction of the big-Oh term in Corollary \ref{1stc}. The size of this fraction might be comparable to $x/\phi(q)$. Indeed, there might exist a sequence $\{q_j\}_{j\in\N}$ of moduli $q$ that correspond to Siegel zeroes larger than $1-o_{j\rightarrow \infty}((\log q_j)^{-1})$. This means that upon using the fact that $L_q(1,\psi)\asymp (1-\beta)\log q$, for $x={q_j}^A$, no matter how large $A$ is, the aforementioned fraction is of size $(1+o_{j\rightarrow \infty}(1)){q_j}^{A}/\phi(q_j)$.
\end{remark}

In \cite[Corollary 1.4]{thza}, Thorner and Zaman proved Theorem \ref{main} with a slightly better error term where $(\log \log x)^{-3/5}$ is replaced by $(\log \log x)^{-1/5}$. The error term of their work is the best to date. Their respective arguments fall in the realm of the classical approaches, as they are based on non-trivial information about the zeroes of the functions $L(\cdot,\chi)$. In this article, we achieve an error term which is almost as strong as theirs by adopting a pretentious approach that avoids a ``heavy" involvement of those zeroes.

We prove Theorem \ref{main} in Section 5. Its proof borrows ideas from the pretentious large sieve that Granville, Harper and Soundararjan developed in \cite{ghal}. It is also inspired by techniques of Koukoulopoulos from his work on bounded multiplicative functions with small partial sums \cite{oldk}. These techniques were deployed in a similar fashion by the author in \cite{sach} where he generalized the work of Koukoulopoulos for classes of divisor-bounded multiplicative functions.

The rest of the text is divided into four sections and the first three serve as preparatory material for the last one where we prove Theorem \ref{main}.

In subsequent work, we intend to extend the recent results of \cite{kousou} and \cite{sach} by applying the core methods of this paper to determine the structure of divisor-bounded multiplicative functions whose partial sums are small on arithmetic progressions.

\subsection*{Notation} Throughout the text, for an integer $n>1$, we denote its smallest prime factor by $P^-(n)$. For $n=1$, we define $P^-(1)=+\infty$. For $m\in\N$, the symbol $\tau_m$ will denote the $m$-fold divisor function given as $\tau_m(n)=\sum_{d_1\cdots d_m=n}1$ for all $n\in\N$. 

\subsection*{Acknowledgements} The author would like to thank his advisor, Dimitris Koukoulopoulos, for his constant guidance and for his suggestions that led to a strengthening of Theorem \ref{main}. He would also like to express his gratitude towards the Stavros Niarchos Foundation for all its generous financial support during the making of this work.

\vspace{1mm}

\section{Preliminary Sieving Results}

In this section we prove two auxiliary sieving results. These are Lemmas \ref{chlem} and \ref{logap}. For the proof of the former, we will use the following lemma.

\begin{lem}\label{lemdb}
Given $q\in\N$, let $\chi$ be a Dirichlet character modulo $q$. For $t\in\R$ and real numbers $x\geq y\geq \max\{q(\abs{t}+1),10\}^{100}$, we have

$$\sum_{\substack{n\leq x\\P^-(n)>y}}\chi(n)n^{it}=\mathds{1}_{\chi=\chi_0}\cdot\frac{x^{1+it}}{1+it}\prod_{p\leq y}\bigg(1-\frac{1}{p}\bigg)+O\bigg(\frac{x^{1-110/\log y}}{\log y}\bigg),$$

\vspace{2mm}

\noindent
where $\chi_0$ is the principal character modulo $q$.
\end{lem}

\begin{proof}
See \cite[Lemma 22.2, p. 224]{dimb}.
\end{proof}

The following theorem of Shiu \cite{shiu} will be necessary for the proofs of Lemmas \ref{chlem} and \ref{logap}. This theorem may be seen as an analogue of the Brun-Titchmarsch inequality for multiplicative functions.

\begin{thm}\label{shiu}
Fix $m\in\N$ and $\eps>0$. Given any choice of $q\in\N, a\in(\Z/q\Z)^*,$ real numbers $x\geq y\geq 1$ with $y/q\geq x^{\eps}$ and a multiplicative function $f$ such that $0\leq f\leq \tau_m$, we have that

$$\sum_{\substack{x-y<n\leq x\\n\equiv a\mod{q}}}f(n)\ll_{m,\eps}\frac{y}{q}\exp\bigg\{\sum_{\substack{p\leq x\\p\nmid q}}\frac{f(p)-1}{p}\bigg\}.$$
\end{thm}

\vspace{1mm}

\begin{lem}\label{chlem}
Given a $q\in\N$, let $\chi$ be a Dirichlet character modulo $q$. For $j\in\N\cup\{0\}$ and any real numbers $x\geq y\geq (10q)^{100}$, we have

$$\sum_{\substack{n\leq x\\P^-(n)>y}}\chi(n)(\log n)^j=\mathds{1}_{\chi=\chi_0}\cdot\bigg(\int_y^x(\log t)^j\dee t\bigg)\prod_{p\leq y}\bigg(1-\frac{1}{p}\bigg)+O\bigg(\frac{(\log x)^jx^{1-\kappa/\log y}}{\log y}\bigg),$$

\noindent
where $\chi_0$ is the principal character modulo $q$ and $\kappa>0$ is an absolute constant.
\end{lem}

\begin{proof}
For any $w\geq y$, Lemma \ref{lemdb} with $t=0$ gives

\begin{equation}\label{chi1}
\sum_{\substack{n\leq w\\P^-(n)>y}}\chi(n)=\mathds{1}_{\chi=\chi_0}\cdot w\prod_{p\leq y}\bigg(1-\frac{1}{p}\bigg)+R_y(w),
\end{equation}

\vspace{2mm}

\noindent
where $R_y(w)\ll(\log y)^{-1}w^{1-\kappa/\log y}$ for some sufficiently small constant $\kappa>0$. Using the Riemann-Stieltjes integral, we have

\begin{align}\label{st}
\sum_{\substack{n\leq x\\P^-(n)>y}}\chi(n)(\log n)^j&=\int_y^x(\log t)^j\dee\Big(\sum_{\substack{n\leq t\\P^-(n)>y}}\chi(n)\Big)\\
&=\mathds{1}_{\chi=\chi_0}\cdot\bigg(\int_y^x(\log t)^j\dee t\bigg)\prod_{p\leq y}\bigg(1-\frac{1}{p}\bigg)+\int_y^x(\log t)^j\dee R_y(t).\nonumber
\end{align}

\vspace{1mm}

\noindent
Now, if $x\geq y^A$ for some large $A>0$, so that $y\leq x^{1-\kappa/\log 2}\leq x^{1-\kappa/\log y}$, then a simple integration by parts implies that

\begin{align*}
\int_y^x(\log t)^j\dee R_y(t)\ll \frac{(\log x)^jx^{1-\kappa/\log y}}{\log y}+\frac{j}{\log y}\int_y^x(\log t)^{j-1}t^{-\kappa/\log y}\dee t.
\end{align*}

\vspace{2mm}

\noindent
But, since

\begin{align*}
j\int_y^x(\log t)^{j-1}t^{-\kappa/\log y}\dee t\leq jx^{1-\kappa/\log y}\int_1^x\frac{(\log t)^{j-1}}{t}\dee t=(\log x)^jx^{1-\kappa/\log y},
\end{align*}

\vspace{2mm}

\noindent
we deduce that

\vspace{-2mm}

\begin{align}\label{est}
\int_y^x(\log t)^j\dee R_y(t)\ll \frac{(\log x)^jx^{1-\kappa/\log y}}{\log y}.
\end{align}

\noindent
We insert (\ref{est}) in (\ref{st}) and complete the proof of the lemma when $x\geq y^A$. It only remains to establish the lemma in the case $y\leq x<y^A$. In this case, $\log x/\log y\asymp 1$, and so applying Theorem \ref{shiu} with $f=\mathds{1}_{P^-(\cdot)>y}$ and $q=a=1$, we have that

$$\bigg|\sum_{\substack{n\leq x\\P^-(n)>y}}\chi(n)(\log n)^j\bigg|\leq (\log x)^j\sum_{\substack{n\leq x\\P^-(n)>y}}1\ll \frac{x(\log x)^j}{\log y}\ll \frac{(\log x)^jx^{1-\kappa/\log y}}{\log y}.$$

\noindent
This means that the lemma does hold in the range $y\leq x<y^A$ as well and this finishes the proof.
\end{proof}

For the proof of Lemma \ref{logap}, we will need Theorem \ref{shiu} as well as an application of the fundamental lemma of sieve theory.

\begin{thm}[The Fundamental Lemma of Sieve Theory]\label{FL}
Let $\CA=\{a_n\}_{n\in \N}$ be a sequence of non-negative real numbers with $\sum_{n\geq 1}a_n<\infty$ and let $\CP$ be a set of primes. Let also $y\geq 1$ be a real number and set $P(y):=\prod_{p\leq y,\ p\in\CP}p$. We define

\vspace{1mm}

\begin{equation*}
S(\CA,\CP,y):=\!\sum_{(n,P(y))=1}\!a_n,
\end{equation*}

\noindent
and for $d\mid P(y)$, we put

\begin{equation*}
A_d:=\sum_{d\mid n}a_n.
\end{equation*}

\vspace{2mm}

\noindent
If there exists a non-negative multiplicative function $v$, some real number $X$, remainder terms $r_d$ and positive constants $\kappa$ and $C$ such that
\vspace{1mm}
\begin{itemize}
    \item $v(p)<p\quad \text{for all}\quad p\mid P(z)$,
    \vspace{2mm}
    \item $\displaystyle{A_d=X\cdot\frac{v(d)}{d}+r_d} \quad \text{for all} \quad d\mid P(z), \quad \text{and}$
    \vspace{2mm}
    \item $\displaystyle{\prod_{\substack{p\in \CP\\y_1<p\leq y_2}}\left(1-\frac{v(p)}{p}\right)^{-1}\leq\left(\frac{\log y_2}{\log y_1}\right)^{\kappa}\left(1+\frac{C}{\log y_1}\right)} \quad \text{for}\quad 2\leq y_1\leq y_2<z,$
\end{itemize}
\vspace{1mm}
then, for every real real number $u\geq 1$, we have that

$$S(\CA,\CP,y)=X\prod_{p\mid P(y)}\left(1-\frac{v(p)}{p}\right)(1+O_{\kappa, C}(u^{-u/2}))+O\Bigg(\sum_{\substack{d\mid P(y)\\d\leq y^u}}|r_d|\Bigg).$$
\end{thm}

\begin{proof}
See \cite[Theorem 18.11, p. 190]{dimb}.
\end{proof}

\vspace{1mm}

\begin{lem}\label{logap}
Let $j$ be a non-negative integer. For any $q\in\N,\, a\in(\Z/q\Z)^*$ and real numbers $x\geq y\geq 2q^2$, there exists an absolute constant $\lambda>0$ such that

\begin{equation}\label{lem2.3}
\sum_{\substack{n\leq x,P^-(n)>y\\n\equiv a \mod{q}}}(\log n)^j=\frac{1}{\phi(q)}\bigg(\int_y^x(\log t)^j\dee t\bigg)\prod_{p\leq y}\bigg(1-\frac{1}{p}\bigg)+O\bigg(\frac{(\log x)^jx^{1-\lambda/\log y}}{\phi(q)\log y}\bigg).
\end{equation}
\end{lem}

\begin{proof}
Let $A>0$ be a sufficiently large real number. We are going to prove the lemma separately in the ranges $x\geq y^A$ and $y\leq x<y^A$.

First, we start with the case $x\geq y^A$. In this case, we will apply Theorem \ref{FL} with $\CP=\{p\leq y, p\nmid q\}$ and $\CA=\{a_n\}_{n\in \N}$, where

\vspace{-1mm}

\begin{eqnarray*}
a_n=\begin{cases}
(\log n)^j,& \text{when} \,\, n\leq x \,\, \text{and} \,\, n\equiv a\mod{q}\\
0,& \text{otherwise}.
\end{cases}
\end{eqnarray*}

\vspace{1mm}

\noindent
Observe that

\begin{equation*}
S(\CA,\CP,y)=\sum_{\substack{n\leq x,P^-(n)>y\\n\equiv a \mod{q}}}(\log n)^j.
\end{equation*}

\vspace{1mm}

\noindent
Consequently, with the chosen sets $\CA$ and $\CP$, an application of Theorem \ref{FL} will provide an asymptotic formula for the sums of interest at the right-hand side.

Now, let $d$ be a positive integer dividing $\prod_{p\leq y,\,p\nmid q}p$. Because of the Chinese Remainder Theorem, the system of linear congruences $n\equiv a\mod{q}, n\equiv 0\mod{d}$ is equivalent to $n\equiv a^*\mod{qd}$ for some $a^*\in(\Z/(qd)\Z)$. Since $\sum_{n\leq x,n\equiv a^*\mod{qd}}1=x/(qd)+O(1)$, partial summation implies that

\begin{align}\label{arel}
\CA_d=\!\sum_{\substack{n\leq x,d\mid n\\n\equiv a\mod{q}}}\!(\log n)^j=\!\sum_{\substack{n\leq x\\n\equiv a^*\mod{qd}}}\!(\log n)^j&=\frac{1}{qd}\bigg\{x(\log x)^j-j\int_1^x(\log t)^{j-1}\dee t\bigg\}+O((\log x)^j)\nonumber\\
&=\frac{1}{qd}\int_1^x(\log t)^j\dee t+O((\log x)^j).
\end{align}

\vspace{2mm}

\noindent
Therefore, following the notation of Theorem \ref{FL}, we have $X=q^{-1}\int_1^x(\log t)^j\dee t$, $v(d)=1$ and $r_d=(\log x)^j$ for all $d\mid \prod_{p\leq y,\,p\nmid q}p$. Then Mertens' theorem implies that we may choose $\kappa=1$ and some large $C>0$ in Theorem \ref{FL}. 

For $u=e\log x/(A\log y)\geq e$, note that

\begin{eqnarray}\label{esti1}
X\prod_{p\leq y\,p\nmid q}\left(1-\frac{v(p)}{p}\right)u^{-u/2}\!\!\!&<&\!\!\!\frac{x(\log x)^j}{q}\prod_{p\leq y}\left(1-\frac{1}{p}\right)\prod_{p\mid q}\left(1-\frac{1}{p}\right)^{-1}\!\exp\Big\{-\frac{e\log x}{2A\log y}\Big\}\nonumber\\
&\leq&\!\!\!\frac{(\log x)^jx^{1-\lambda/\log y}}{\phi(q)}\prod_{p\leq y}\left(1-\frac{1}{p}\right)\nonumber\\
&\ll&\!\!\!\frac{(\log x)^jx^{1-\lambda/\log y}}{\phi(q)\log y}
\end{eqnarray}

\vspace{1mm}

\noindent
for some $\lambda\in(0,e(2A)^{-1}].$ Furthermore, with $P(y)=\prod_{p\leq y,\, p\nmid q}p$ and the choice of $u$ that we made, we have that

\begin{eqnarray}\label{esti2}
\sum_{\substack{d\mid P(y)\\d\leq y^u}}|r_d|\leq \sum_{d\leq y^u}r_d\!\!\!&\leq&\!\!\!x^{e/A}(\log x)^j\ll x^{2e/A}(\log x)^{j-1}\nonumber\\
&<&\!\!\!\frac{(\log x)^jx^{(2e+1)/A}}{\phi(q)\log y}\leq\frac{(\log x)^jx^{1-\lambda/\log y}}{\phi(q)\log y},
\end{eqnarray}

\vspace{1mm}

\noindent
where we passed from the first line to the second by using the inequality $x^{1/A}\geq y\geq2q>\phi(q)$. The last step follows from the fact that $A$ is sufficiently large.

Combination of (\ref{esti1}) and (\ref{esti2}) with Theorem \ref{FL} yields

\vspace{-1mm}

\begin{align}\label{int}
\sum_{\substack{n\leq x,P^-(n)>y\\n\equiv a \mod{q}}}(\log n)^j=\frac{1}{\phi(q)}\bigg(\int_1^x(\log t)^j\dee t\bigg)\prod_{p\leq y}\bigg(1-\frac{1}{p}\bigg)+O\bigg(\frac{(\log x)^jx^{1-\lambda/\log y}}{\phi(q)\log y}\bigg).
\end{align}

\vspace{1mm}

\noindent
But, $\prod_{p\leq y}(1-1/p)\ll (\log y)^{-1}$, and

\vspace{-1mm}

\begin{align*}
\int_1^y(\log t)^j\dee t<y(\log x)^j\ll (\log x)^jx^{1-\lambda/\log y},
\end{align*}

\vspace{1mm}

\noindent
since $x\geq y^A$ with $A$ sufficiently large. Hence, relation (\ref{int}) implies (\ref{lem2.3}) in the range $x\geq y^A$. 

When $y\leq x<y^A$, a similar argument as the one that we developed at the end of the proof of Lemma \ref{chlem} shows that

\vspace{-1mm}

$$\sum_{\substack{n\leq x,P^-(n)>y\\n\equiv a \mod{q}}}(\log n)^j\ll\frac{(\log x)^jx^{1-\lambda/\log y}}{\phi(q)\log y}$$

\vspace{1mm}

\noindent
and the proof of the lemma is complete.
\end{proof}

\vspace{1mm}

\section{Existing Bounds for Sifted $L$-Dirichlet Series}

Let us consider a real number $y\geq 1$ and a positive integer $q$. For $s\in\C$ with $\Re(s)>1$, we define the $y$-rough Dirichlet series of a Dirichlet character $\chi$ modulo $q$ as

\begin{equation}\label{sifL}
L_y(s,\chi):=\sum_{\substack{n\geq 1\\P^-(n)>y}}\frac{\chi(n)}{n^s}=L(s,\chi)\prod_{p\leq y}(1-\chi(p)/p^s).
\end{equation}

\noindent
In (\ref{sifL}), the series is absolutely convergent and its rightmost side implies the meromorphical continuation of $L_y(\cdot,\chi)$ on the whole complex plane (with one pole at $1$ only in the case of the principal character $\chi_0$).

The present section constitutes a collection of some existing upper and lower bounds for the values of a $y$-rough Dirichlet series. We will use all these bounds in the proof of Theorem \ref{main}.
The theorems listed below are stated without their proofs, as they are already part of the literature. We only provide a reference for each one of them. To simplify their statements, we introduce the notation

$$V_t:=\exp\{100(\log (3+|t|))^{2/3}(\log\log(3+|t|))^{1/3}\}\quad (t\in\R).$$

\vspace{3mm}

We start with the following theorem \cite[Lemma 4.1]{ppnt} concerning an upper bound for the derivatives of $L(\cdot,\chi)$.

\begin{thm}\label{ub}
Let $q$ be a positive integer and $\chi$ be a non-principal character modulo $q$. Let also $j\in\N$ and $s=\sigma+it$ with $\sigma>1$ and $t\in\R$. For $y\geq qV_t$, we have that

$$|L_y^{(j)}(s,\chi)|\ll j!(C\log y)^j,$$

\vspace{2mm}

\noindent
where $C>0$ is an absolute constant.
\end{thm}

\vspace{1mm}

The next result \cite[Lemma 4.2]{ppnt} equips us with lower bounds for $L_y(s,\chi)$.

\begin{thm}\label{lb}
Fix a positive integer $q$ and let $\chi$ be a character modulo $q$. Let $s=\sigma+it$ with $\sigma>1$ and $t\in\R$ and consider the real number $y\geq qV_t$.

\begin{enumerate}[(a)]
\item If $\chi$ is not real, then $\abs{L_y(s,\chi)}\gg 1$.
\item If $\chi$ is real and non-principal, then $\abs{L_y(s,\chi)}\gg L_y(1,\chi).$
\end{enumerate}
\end{thm}

\vspace{1mm}

We close this section with a theorem \cite[Lemma 27.5, p. 291]{dimb} dealing with the size of $L_q(\sigma,\chi)$ for $\sigma\geq 1$ when $\chi\neq\psi$, where $\psi$ is defined as in the statement of Theorem \ref{main}.

\begin{thm}\label{ldnS}
Let $q$ be a positive integer and let $\CR_q$ be the set of real, non-principal characters modulo $q$. If we take a character $\psi$ such that $L_q(1,\psi)=\min_{\chi\in\CR_q}L_q(1,\chi)$ and $\CC_q:=\{\chi\mod{q}:\chi\neq \chi_0, \psi\}$, then $|L_y(\sigma,\chi)|\asymp 1$ for all $\chi\in\CC_q, y\geq q$ and $\sigma\geq 1$.
\end{thm}

\vspace{1mm}

\section{A Mean Value Theorem}

The goal of this section is to establish a mean value theorem for the derivatives of the Dirichlet series 

$$\sum_{\substack{n\equiv a\mod{q}\\P^-(n)>y}}\frac{\Lambda(n)}{n^s}.$$

\vspace{2mm}

\noindent
Such a result will be essential at the end of the proof of Theorem \ref{main}. To prove it, we are making use of the following lemma which is due to Montgomery \cite[Theorem 3, p. 131]{mon}.

\begin{lem}\label{mon}
Let $A(s)=\sum_{n\geq 1}a_nn^{-s}$ and $B(s)=\sum_{n\geq 1}b_nn^{-s}$ be two Dirichlet series which converge for $\Re(s)>1.$ If $|a_n|\leq b_n$ for all $n\in \N$, then

\begin{equation*}
\int_{-T}^T|A(\sigma+it)|^2\dee t\leq 3\int_{-T}^T|B(\sigma+it)|^2\dee t,
\end{equation*}

\vspace{3mm}
	
\noindent
for any $\sigma>1$ and any $T\geq 0$. 
\end{lem}

\begin{lem}\label{mvt}
Let us consider an integer $q\geq 1$ and a real number $y\geq 16q^2$. For $j\in\N\cup\{0\}, T\geq 1, \sigma\in(1,2)$ and $a\in(\Z/q\Z)^*$, we have
	
$$\int_{|t|>T}\bigg|\sum_{\substack{n\equiv a\mod{q}\\P^-(n)>y}}\frac{\Lambda(n)(\log n)^j}{n^{\sigma+it}}\bigg|^2\frac{\dee t}{t^2}\ll\frac{16^{\,j}(2j)!}{\phi(q)^2(\sigma-1)^{2j+1}T}.$$
\end{lem}

\begin{proof}
First, for $k\in\N\cup\{0\}$, because of Lemma \ref{mon}, we observe that

\begin{align*}
\int_{k-1/2}^{k+1/2}\bigg|\sum_{\substack{n\equiv a\mod{q}\\P^-(n)>y}}\frac{\Lambda(n)(\log n)^j}{n^{\sigma+it}}\bigg|^2\dee t&=\int_{-1/2}^{1/2}\bigg|\sum_{\substack{n\equiv a\mod{q}\\P^-(n)>y}}\frac{\Lambda(n)(\log n)^jn^{-ik}}{n^{\sigma+it}}\bigg|^2\dee t\\
&\leq 3\int_{-1/2}^{1/2}\bigg|\sum_{\substack{n\equiv a\mod{q}\\P^-(n)>y}}\frac{\Lambda(n)(\log n)^j}{n^{\sigma+it}}\bigg|^2\dee t.\nonumber
\end{align*}

\noindent
Consequently,

\begin{align}\label{lines}
\int_{|t|>T}\bigg|\sum_{\substack{n\equiv a\mod{q}\\P^-(n)>y}}\frac{\Lambda(n)(\log n)^j}{n^{\sigma+it}}\bigg|^2\frac{\dee t}{t^2}&\leq \sum_{|k|>T-1/2}\int_{k-1/2}^{k+1/2}\bigg|\sum_{\substack{n\equiv a\mod{q}\\P^-(n)>y}}\frac{\Lambda(n)(\log n)^j}{n^{\sigma+it}}\bigg|^2\frac{\dee t}{t^2}\nonumber\\
&\leq 4\sum_{|k|>T/2}\frac{1}{k^2}\int_{k-1/2}^{k+1/2}\bigg|\sum_{\substack{n\equiv a\mod{q}\\P^-(n)>y}}\frac{\Lambda(n)(\log n)^j}{n^{\sigma+it}}\bigg|^2\dee t\nonumber\\
&\ll\bigg(\sum_{k>T/2}\frac{1}{k^2}\bigg)\cdot \int_{-1/2}^{1/2}\bigg|\sum_{\substack{n\equiv a\mod{q}\\P^-(n)>y}}\frac{\Lambda(n)(\log n)^j}{n^{\sigma+it}}\bigg|^2\dee t\nonumber\\
&\ll\frac{1}{T}\int_{-1/2}^{1/2}\bigg|\sum_{\substack{n\equiv a\mod{q}\\P^-(n)>y}}\frac{\Lambda(n)(\log n)^j}{n^{\sigma+it}}\bigg|^2\dee t.
\end{align}

Now, we focus on estimating the integral at the last line of (\ref{lines}). We will do this by adopting the rather standard technique which is used for proving similar mean value theorems. We consider the function $\Phi:\R\rightarrow \R$ given by the formula $\Phi(x)=(2\pi\sin(x/4))^2x^{-2}$ for all $x\in\R^*$ and $\Phi(0)=\pi^2/4$. Notice that $\Phi(x)\geq 1$ on $[-1/2,1/2]$. So, if $\hat{\Phi}$ is the Fourier transform of $\Phi,$ then

\begin{align}\label{afw}
\int_{-1/2}^{1/2}\bigg|\sum_{\substack{n\equiv a\mod{q}\\P^-(n)>y}}\frac{\Lambda(n)(\log n)^j}{n^{\sigma+it}}\bigg|^2\dee t&\leq \int_{\R}\bigg|\sum_{\substack{n\equiv a\mod{q}\\P^-(n)>y}}\frac{\Lambda(n)(\log n)^j}{n^{\sigma+it}}\bigg|^2\Phi(t)\dee t\nonumber\\
&=\sum_{\substack{m\equiv a\mod{q}\\n\equiv a\mod{q}\\P^-(m),P^-(n)>y}}\frac{\Lambda(n)\Lambda(n)(\log m)^j(\log n)^j}{m^{\sigma}n^{\sigma}}\hat{\Phi}(\log(m/n)),
\end{align}

\noindent
where we arrived at the last line by expanding the square and by interchanging the order of summation and integration. The Fourier transform $\hat{\Phi}$ is an even function, because $\Phi$ is also even. Therefore, we may bound the last line of (\ref{afw}) by twice the sum

\begin{equation}\label{sum}
\sum_{\substack{m\equiv a\mod{q}\\n\leq m, n\equiv a\mod{q}\\P^-(m),P^-(n)>y}}\frac{\Lambda(n)\Lambda(n)(\log m)^j(\log n)^j}{m^{\sigma}n^{\sigma}}\hat{\Phi}(\log(m/n)).
\end{equation} 

\noindent
The Fourier transform $\hat{\Phi}$ is continuous and compactly supported on $[-1/2,1/2]$. Moreover, for $n\leq m$, we have that $|n-m|\leq m\log(m/n)$. Therefore, the sum (\ref{sum}) is smaller than or equal to 

\begin{align*}
\sum_{\substack{m\equiv a\mod{q}\\P^-(m)>y}}\frac{\Lambda(m)(\log m)^j}{m^{\sigma}}&\sum_{\substack{n\leq m, |n-m|\leq m/2\\n\equiv a\mod{q}}}\frac{\Lambda(n)(\log n)^j}{n^{\sigma}}\\
&\ll\sum_{\substack{m\equiv a\mod{q}\\P^-(m)>y}}\frac{\Lambda(m)(\log m)^{2j}}{m^{2\sigma}}\sum_{\substack{m/2\leq n\leq m\\n\equiv a\mod{q}}}\Lambda(n)\\
&\ll\frac{1}{\phi(q)}\sum_{\substack{m\equiv a\mod{q}\\P^-(m)>y}}\frac{\Lambda(m)(\log m)^{2j}}{m^{2\sigma-1}}.
\end{align*}

\noindent
For the last estimate we made use of the Brun-Titchmarsch theorem. Its application was allowed, because $\Lambda(1)=0$, which means that the condition $P^-(m)>y$ implies that $m>y\geq 16q^2$. According to all the above, relation (\ref{afw}) becomes

\begin{equation}\label{1stbe}
\int_{-1/2}^{1/2}\bigg|\sum_{\substack{n\equiv a\mod{q}\\P^-(n)>y}}\frac{\Lambda(n)(\log n)^j}{n^{\sigma+it}}\bigg|^2\dee t\ll\frac{1}{\phi(q)}\sum_{\substack{m\equiv a\mod{q}\\P^-(m)>y}}\frac{\Lambda(m)(\log m)^{2j}}{m^{2\sigma-1}}.
\end{equation}

We continue by bounding the sum of the right-hand side of (\ref{1stbe}). By decomposing this sum in dyadic intervals, we get that

\begin{align*}
\sum_{\substack{m\equiv a\mod{q}\\P^-(m)>y}}\frac{\Lambda(m)(\log m)^{2j}}{m^{2\sigma-1}}&\leq (\log 4)^j\sum_{r>\log y/\log 2}\frac{r^{2j}}{2^{r(2\sigma-1)}}\sum_{\substack{2^r\leq m<2^{r+1}\\m\equiv a\mod{q}}}\Lambda(m)\\
&\ll\frac{(\log 4)^j}{\phi(q)}\sum_{r\geq 0}\frac{r^{2j}}{4^{r(\sigma-1)}},
\end{align*}

\noindent
where we applied the Brun-Titchmarsch theorem for the last step, since $2^r>y\geq 16q^2$. Using the $2j$-th derivative of the geometric series, it follows that

\begin{equation*}
\sum_{r\geq 0}\frac{r^{2j}}{4^{r(\sigma-1)}}\leq \sum_{r\geq 0}\frac{(r+2j)!}{4^{r(\sigma-1)}r!}=\frac{(2j)!}{(1-4^{1-\sigma})^{2j+1}}\ll \frac{(2j)!(2/(\log 2))^{2j}}{(\sigma-1)^{2j+1}},
\end{equation*}

\vspace{2mm}

\noindent
as $1-4^{1-\sigma}\geq \frac{(\sigma-1)\log 2}{2}$ for $\sigma\in(1,2)$ by the mean value theorem. Thus,

\begin{equation}\label{ests}
\sum_{\substack{m\equiv a\mod{q}\\P^-(m)>y}}\frac{\Lambda(m)(\log m)^{2j}}{m^{2\sigma-1}}\ll\frac{(8/\log 2)^j(2j)!}{\phi(q)(\sigma-1)^{2j+1}}<\frac{16^{\,j}(2j)!}{\phi(q)(\sigma-1)^{2j+1}}.
\end{equation}

We now combine (\ref{lines}), (\ref{1stbe}) and (\ref{ests}) and conclude the proof of the lemma.  
\end{proof}

\vspace{1mm}

\section{Proof of Theorem \ref{main}}

\subsection{A useful lemma} Before embarking on the proof of Theorem \ref{main}, which is the main objective of this section, we first state a lemma that will be useful. 

\begin{lem}\label{lasl}
Let $k\in\N$, $D$ be an open set of $\C$, $s\in D$ and $F:D\rightarrow\C$ be a function which is differentiable $k$ times at $s$. We further assume that $F(s)\neq 0$ and we set

$$K=\max_{1\leq j\leq k}\bigg\{\frac{1}{j!}\bigg|\frac{F^{(j)}}{F}(s)\bigg|\bigg\}^{1/j}\quad and\quad N=\max_{1\leq j\leq k}\bigg\{\frac{1}{j!}\bigg|\Big(\frac{F^{'}}{F}\Big)^{(j-1)}(s)\bigg|\bigg\}^{1/j}.$$

\vspace{2mm}

\noindent
Then $K/2\leq N\leq 2K$.
\end{lem}

\begin{proof}
See \cite[Lemma 9.1]{oldk}.
\end{proof}

\subsection{Proof of Theorem \ref{main}:} We may assume that $q\geq 2$ and that $x\geq q^B$ for some sufficiently large real number $B>0$. Indeed, for $q=1$, we have the prime number theorem with the error term provided by Korobov and Vinogradov and if $q^2\leq x<q^B$, then $\log x/\log q\asymp 1$ and the theorem follows from a trivial application of the Brun-Titchmarsch inequality. Now, let $k$ be a positive integer and set $y=(10q)^{100}V_T$ with $T=\exp\{2L(\log x)^{3/5}(\log\log x)^{2/5}\}$ for a large constant $L>0$. 

Since $\log =\Lambda\ast 1$, we have that

\begin{align}\label{1'}
\sum_{\substack{n\leq x\\n\equiv a \mod{q}\\P^-(n)>y}}\log n=\sum_{\substack{m\leq \sqrt{x}\\P^-(m)>y}}\sum_{\substack{\ell\leq \frac{x}{m}\\\ell\equiv a\overline{m}\mod{q} \\P^-(\ell)>y}}\Lambda(\ell)+\sum_{\substack{\ell\leq \sqrt{x}\\P^-(n)>y}}\Lambda(\ell)\sum_{\substack{\sqrt{x}<m\leq \frac{x}{\ell}\\m\equiv a\overline{\ell}\mod{q}\\P^-(m)>y}}1.
\end{align}

\vspace{2mm}

\noindent
In the rest of the proof, the character $\psi$ is defined as in the statement of Theorem \ref{main}. For $\chi\in\{\chi_0,\psi\}$, it is true that $\log\cdot\chi=(\Lambda\cdot\chi)\ast\chi$. Therefore, we similarly have that

\begin{align}\label{2'}
\sum_{\substack{n\leq x\\P^-(n)>y}}\chi(n)\log n=\sum_{\substack{m\leq \sqrt{x}\\P^-(m)>y}}\chi(m)\sum_{\substack{\ell\leq \frac{x}{m}\\P^-(\ell)>y}}\Lambda(\ell)\chi(\ell)+\sum_{\substack{\ell\leq \sqrt{x}\\P^-(n)>y}}\Lambda(\ell)\chi(\ell)\sum_{\substack{\sqrt{x}<m\leq \frac{x}{\ell}\\P^-(m)>y}}\chi(m).
\end{align}

\vspace{2mm}

\noindent
We now use Lemma \ref{chlem} for the sums $\sum_{n\leq x,\,P^-(n)>y}\chi_0(n)\log n$ and $\sum_{n\leq x,\,P^-(n)>y}\psi(n)\log n$. We also use Lemma \ref{logap} for the sum $\sum_{\substack{n\leq x,\,P^-(n)>y\\n\equiv a \mod{q}}}\log n$. Combination of the obtained formulas yields 

\begin{align}\label{last}
\sum_{\substack{n\leq x\\n\equiv a \mod{q}\\P^-(n)>y}}\log n-\frac{1}{\phi(q)}\sum_{\substack{n\leq x\\P^-(n)>y}}\chi_0(n)\log n-\frac{\psi(a)}{\phi(q)}\sum_{\substack{n\leq x\\P^-(n)>y}}\psi(n)\log n\ll\frac{x^{1-c_1/\log y}}{\phi(q)},
\end{align}

\vspace{2mm}

\noindent
where $c_1=\min\{\kappa,\lambda\}/2$. Now, for simplicity, we put 

\begin{align*}
\Delta(u,y;q,b):=\sum_{\substack{\ell\leq u\\\ell\equiv b\mod{q}\\P^-(\ell)>y}}\Lambda(\ell)-\frac{1}{\phi(q)}\sum_{\substack{\ell\leq u\\P^-(\ell)>y}}\Lambda(\ell)\chi_0(\ell)-\frac{\psi(b)}{\phi(q)}\sum_{\substack{\ell\leq u\\P^-(\ell)>y}}\Lambda(\ell)\psi(\ell)
\end{align*}

\noindent
and

\begin{align*}
\Delta^*(u,y;q,b):=\sum_{\substack{m\leq u\\m\equiv b\mod{q}\\P^-(m)>y}}1-\frac{1}{\phi(q)}\sum_{\substack{m\leq u\\P^-(m)>y}}\chi_0(m)-\frac{\psi(b)}{\phi(q)}\sum_{\substack{m\leq u\\P^-(m)>y}}\psi(m).
\end{align*}

\vspace{2mm}

\noindent
We take (\ref{2'}) once with $\chi=\chi_0$ and once with $\chi=\psi$ and then we add the two relations term by term. Then we subtract the resulting relation from (\ref{1'}). This leads to

\begin{align*}
\sum_{\substack{m\leq\sqrt{x}\\P^-(m)>y}}\Delta\Big(\frac{x}{m},y;q,a\overline{m}\Big)&=-\sum_{\substack{\ell\leq\sqrt{x}\\P^-(\ell)>y}}\Lambda(\ell)\Delta^*\Big(\frac{x}{\ell},y;q,a\overline{\ell}\Big)\\
&+\sum_{\substack{\ell\leq\sqrt{x}\\P^-(\ell)>y}}\Lambda(\ell)\Delta^*(\sqrt{x},y;q,a\overline{\ell})+O\bigg(\frac{x^{1-c_1/\log y}}{\phi(q)}\bigg),
\end{align*}

\vspace{1mm}

\noindent
where the big-Oh term comes from the contribution of the left-hand side of (\ref{last}). Since $x/\ell\geq\sqrt{x}\geq y$ for $\ell\leq \sqrt{x}$, we can apply Lemmas \ref{chlem} and \ref{logap} with $j=0$ to bound the three sums in the definitions of $\Delta^*(x/\ell,y;q,a\overline{\ell})$ and $\Delta^*(\sqrt{x},y;q,a\overline{\ell})$. Doing so yields that

\begin{align*}
\sum_{\substack{\ell\leq\sqrt{x}\\P^-(\ell)>y}}\Lambda(\ell)\Big\{\Delta^*\Big(\frac{x}{\ell},y;q,a\overline{\ell}\Big)-\Delta^*(\sqrt{x},y;q,a\overline{\ell})\Big\}&\ll\frac{x^{1-c_1/\log y}}{\phi(q)\log y}\sum_{\substack{\ell\leq\sqrt{x}\\P^-(\ell)>y}}\frac{\Lambda(\ell)}{{\ell}^{1-c_1/\log y}}\\
&\ll\frac{x^{1-c_1/(2\log y)}}{\phi(q)\log y}\sum_{\ell\leq\sqrt{x}}\frac{\Lambda(\ell)}{\ell}\\
&\ll \frac{x^{1-c_1/(2\log y)}\log x}{\phi(q)\log y}\ll \frac{x^{1-c_1/(3\log y)}}{\phi(q)}.
\end{align*}

\noindent
Hence,

\begin{align}\label{rec}
\Delta(x,y;q,a)=-\sum_{\substack{1<m\leq\sqrt{x}\\P^-(m)>y}}\Delta\Big(\frac{x}{m},y;q,a\overline{m}\Big)+O\bigg(\frac{x^{1-c_2/\log y}}{\phi(q)}\bigg),
\end{align}

\noindent
where $c_2=c_1/3$.

Now, we set $D=x^{1-\delta/\log y}$ for some sufficiently small $\delta\in(0,c_2)$ which will be chosen later. When $m\leq \sqrt{x}$ and $x\geq q^B$ with a suitably large exponent $B$, we have that $D/m\geq 2q\sqrt{x/m}$. Thus, for $t\in[x-D,x]$, an application of the Brun-Titchmarsch theorem gives

$$\sum_{\substack{\frac{t}{m}<\ell\leq\frac{x}{m}\\\ell\equiv a\overline{m}\mod{q}\\ P^-(\ell)>y}}\Lambda(\ell)\leq\sum_{\substack{\frac{x-D}{m}<\ell\leq\frac{x}{m}\\\ell\equiv a\overline{m}\mod{q}}}\Lambda(\ell)\ll\frac{D}{m\phi(q)}.$$

\noindent
Similarly,

$$\sum_{\substack{\frac{t}{m}<\ell\leq\frac{x}{m}\\P^-(\ell)>y}}\Lambda(\ell)\chi_0(\ell)\ll\frac{D}{m}\quad\,\,\,\,\text{and}\quad\sum_{\substack{\frac{t}{m}<\ell\leq\frac{x}{m}\\P^-(\ell)>y}}\Lambda(\ell)\psi(\ell)\ll\frac{D}{m},$$

\vspace{2mm}

\noindent
by bounding the characters trivially before making use of the Brun-Titchmarsh theorem. With these estimates, we deduce that

\begin{align}\label{trick}
\sum_{\substack{1<m\leq\sqrt{x}\\P^-(m)>y}}\Delta\Big(\frac{x}{m},y;q,a\overline{m}\Big)&=\frac{1}{D}\sum_{\substack{1<m\leq\sqrt{x}\\P^-(m)>y}}\int_{x-D}^x\Delta\Big(\frac{t}{m},y;q,a\overline{m}\Big)\dee t\\
&+O\bigg(\frac{D}{\phi(q)}\sum_{\substack{m\leq\sqrt{x}\\P^-(m)>y}}\frac{1}{m}\bigg).\nonumber
\end{align}

\noindent
But,

$$\sum_{\substack{m\leq\sqrt{x}\\P^-(m)>y}}\frac{1}{m}\leq \prod_{y<p\leq\sqrt{x}}\bigg(1-\frac{1}{p}\bigg)^{-1}\ll \frac{\log x}{\log y},$$

\noindent
and so the inequality $x^{\delta/(2\log y)}\gg_{\delta}\log x/\log y$ and combination of (\ref{rec}) and (\ref{trick}) lead to

\vspace{-3mm}

\begin{align}\label{recres}
\Delta(x,y;q,a)=-\frac{1}{D}\sum_{\substack{1<m\leq\sqrt{x}\\P^-(m)>y}}\int_{x-D}^x\Delta\Big(\frac{t}{m},y;q,a\overline{m}\Big)\dee t+O_{\delta}\bigg(\frac{x^{1-\delta/(2\log y)}}{\phi(q)}\bigg).
\end{align}

\noindent 

By referring to the orthogonality of the characters modulo $q$, it follows that

\vspace{-1mm}

$$\Delta\Big(\frac{t}{m},y;q,a\overline{m}\Big)=\frac{1}{\phi(q)}\sum_{\substack{\chi\mod{q}\\\chi\neq\chi_0,\psi}}\overline{\chi}(a)\chi(m)\sum_{\substack{\ell\leq\frac{t}{m}\\P^-(\ell)>y}}\Lambda(\ell)\chi(\ell).$$

\vspace{-4mm}

\noindent
Therefore,

\vspace{-3mm}

\begin{align*}
&\frac{1}{D}\bigg|\sum_{\substack{1<m\leq\sqrt{x}\\P^-(m)>y}}\int_{x-D}^x\Delta\Big(\frac{t}{m},y;q,a\overline{m}\Big)\dee t\bigg|\\
&\leq \frac{1}{D\phi(q)}\sum_{\substack{\chi\mod{q}\\\chi\neq\chi_0,\psi}}\bigg|\sum_{\substack{1<m\leq\sqrt{x}\\P^-(m)>y}}\chi(m)\int_{x-D}^x\Big(\sum_{\substack{\ell\leq \frac{t}{m}\\P^-(\ell)>y}}\Lambda(\ell)\chi(\ell)\Big)\dee t\bigg|\\
&=\frac{1}{D\phi(q)}\sum_{\substack{\chi\mod{q}\\\chi\neq\chi_0,\psi}}\bigg|\sum_{\substack{1<m\leq\sqrt{x}\\P^-(m)>y}}\chi(m)m\int_{\frac{x-D}{m}}^{\frac{x}{m}}\Big(\sum_{\substack{\ell\leq t\\P^-(\ell)>y}}\Lambda(\ell)\chi(\ell)\Big)\dee t\bigg|\\
&=\frac{1}{D\phi(q)}\sum_{\substack{\chi\mod{q}\\\chi\neq\chi_0,\psi}}\bigg|\sum_{\substack{1<m\leq\sqrt{x}\\P^-(m)>y}}\chi(m)m\int_{\frac{x-D}{\sqrt{x}}}^{\frac{x}{y}}\mathds{1}_{\big(\frac{x-D}{m},\frac{x}{m}\big]}(t)\cdot\Big(\sum_{\substack{\ell\leq t\\P^-(\ell)>y}}\Lambda(\ell)\chi(\ell)\Big)\dee t\bigg|\\
&=\frac{1}{D\phi(q)}\sum_{\substack{\chi\mod{q}\\\chi\neq\chi_0, \psi}}\bigg|\int_{\frac{x-D}{\sqrt{x}}}^{\frac{x}{y}}\Big(\sum_{\substack{\ell\leq t\\P^-(\ell)>y}}\Lambda(\ell)\chi(\ell)\Big)\Big(\sum_{\substack{\frac{x-D}{t}<m\leq\frac{x}{t}\\P^-(m)>y}}\chi(m)m\Big)\dee t\bigg|.
\end{align*}

\noindent
We move the absolute value inside the integral and then use the Cauchy-Schwarz inequality twice to obtain

\vspace{-3mm}

\begin{align}\label{dcs}
&\frac{1}{D}\bigg|\sum_{\substack{1<m\leq\sqrt{x}\\P^-(m)>y}}\int_{x-D}^x\Delta\Big(\frac{t}{m},y;q,a\overline{m}\Big)\dee t\bigg|\\
&\leq\frac{1}{D\phi(q)}\int_{\frac{x-D}{\sqrt{x}}}^{\frac{x}{y}}\bigg(\sum_{\substack{\chi\mod{q}\\\chi\neq \chi_0, \psi}}\Big|\sum_{\substack{\ell\leq t\\P^-(\ell)>y}}\Lambda(\ell)\chi(\ell)\Big|^2\bigg)^{1/2}\bigg(\sum_{\substack{\chi\mod{q}\\\chi\neq \chi_0, \psi}}\Big|\sum_{\substack{\frac{x-D}{t}<m\leq\frac{x}{t}\\P^-(m)>y}}\chi(m)m\Big|^2\bigg)^{1/2}\dee t\nonumber\\
&\leq\frac{1}{D\phi(q)}\bigg(\sum_{\substack{\chi\mod{q}\\\chi\neq \chi_0, \psi}}\int_{\frac{x-D}{\sqrt{x}}}^{\frac{x}{y}}\Big|\sum_{\substack{\ell\leq t\\P^-(\ell)>y}}\Lambda(\ell)\chi(\ell)\Big|^2\frac{\dee t}{t^3}\bigg)^{1/2}\nonumber\\
&\times \bigg(\int_{\frac{x-D}{\sqrt{x}}}^{\frac{x}{y}}\sum_{\chi\mod{q}}\Big|\sum_{\substack{\frac{x-D}{t}<m\leq\frac{x}{t}\\P^-(m)>y}}\chi(m)m\Big|^2t^3\dee t\bigg)^{1/2}.\nonumber
\end{align}

\noindent
However,

\vspace{-1mm}

\begin{align*}
\sum_{\chi\mod{q}}\Big|\sum_{\substack{\frac{x-D}{t}<m\leq\frac{x}{t}\\P^-(m)>y}}\chi(m)m\Big|^2&=\phi(q)\sum_{b\in(\Z/q\Z)^*}\bigg(\sum_{\substack{\frac{x-D}{t}<m\leq \frac{x}{t}\\m\equiv b\mod{q}\\P^-(m)>y}}m\bigg)^2\\
&\leq \frac{\phi(q)x^2}{t^2}\sum_{b\in(\Z/q\Z)^*}\bigg(\sum_{\substack{\frac{x}{2t}<m\leq\frac{x}{t}\\m\equiv b\mod{q}\\P^-(m)>y}}1\bigg)^2\ll\frac{x^4}{t^4(\log y)^2},
\end{align*}

\vspace{1mm}

\noindent
where the final step is a consequence of Theorem \ref{shiu} which was applied separately to each enlarged sum of the last line. So,

\vspace{-1mm}

$$\int_{\frac{x-D}{\sqrt{x}}}^{\frac{x}{y}}\sum_{\chi\mod{q}}\Big|\sum_{\substack{\frac{x-D}{t}<m\leq\frac{x}{t}\\P^-(m)>y}}\chi(m)m\Big|^2t^3\dee t\ll\frac{x^4}{(\log y)^2}\int_1^x\frac{\dee t}{t}=\frac{x^4\log x}{(\log y)^2}.$$

\vspace{1mm}

\noindent
This means that (\ref{dcs}) leads to the estimate

\vspace{-1mm}

\begin{align}\label{bP}
\frac{1}{D}\sum_{\substack{1<m\leq\sqrt{x}\\P^-(m)>y}}\int_{x-D}^x\Delta\Big(\frac{t}{m},y;q,a\overline{m}\Big)\dee t&\ll \frac{x^{1+\delta/\log y}\sqrt{\log x}}{\phi(q)\log y}\\
&\times\bigg(\sum_{\substack{\chi\mod{q}\\\chi\neq \chi_0, \psi}}\int_{\frac{\sqrt{x}}{2}}^x\Big|\sum_{\substack{\ell\leq t\\P^-(\ell)>y}}\Lambda(\ell)\chi(\ell)\Big|^2\frac{\dee t}{t^3}\bigg)^{1/2}.\nonumber
\end{align}

\vspace{1mm}

Now, we continue by establishing a bound for the sum of integrals at the second line of (\ref{bP}). For $t\in[\sqrt{x}/2,x]$, by partial summation and Chebyshev's estimates, we have that

\vspace{-2mm}

\begin{align*}
\sum_{\substack{\ell\leq t\\P^-(\ell)>y}}\Lambda(\ell)\chi(\ell)&=O(\sqrt{t})+\int_{\sqrt{t}}^t(\log u)^{-k}\dee \Big(\sum_{\substack{\ell\leq u\\P^-(\ell)>y}}\Lambda(\ell)\chi(\ell)(\log \ell)^k\Big)\\
&=O(\sqrt{t})+(\log t)^{-k}\sum_{\substack{\ell\leq t\\P^-(\ell)>y}}\Lambda(\ell)\chi(\ell)(\log \ell)^k\\
&+k\int_{\sqrt{t}}^t\sum_{\substack{\ell\leq u\\P^-(\ell)>y}}\Lambda(\ell)\chi(\ell)(\log \ell)^k\frac{\dee u}{u(\log u)^{k+1}}\\
&\ll\sqrt{t}+M^k(\log x)^{-k}\Big|\sum_{\substack{\ell\leq t\\P^-(\ell)>y}}\Lambda(\ell)\chi(\ell)(\log \ell)^k\Big|\\
&+M^k(\log x)^{-(k+1)}\int_{\sqrt{t}}^t\Big|\sum_{\substack{\ell\leq u\\P^-(\ell)>y}}\Lambda(\ell)\chi(\ell)(\log \ell)^k\Big|\frac{\dee u}{u},
\end{align*}

\noindent
where $M$ is some sufficiently large positive constant. Because of the basic inequality $3(\alpha^2+\beta^2+\gamma^2)\geq (\alpha+\beta+\gamma)^2$ for all $\alpha, \beta, \gamma\in\R$, we deduce that

\vspace{-2mm}

\begin{align*}
\Big|\sum_{\substack{\ell\leq t\\P^-(\ell)>y}}\Lambda(\ell)\chi(\ell)\Big|^2&\ll t+M^{2k}(\log x)^{-2k}\Big|\sum_{\substack{\ell\leq t\\P^-(\ell)>y}}\Lambda(\ell)\chi(\ell)(\log \ell)^k\Big|^2\\
&+M^{2k}(\log x)^{-2(k+1)}\bigg(\int_{\sqrt{t}}^t\Big|\sum_{\substack{\ell\leq u\\P^-(\ell)>y}}\Lambda(\ell)\chi(\ell)(\log \ell)^k\Big|\frac{\dee u}{u}\bigg)^2.
\end{align*}

\noindent
But, upon noticing that

\vspace{-1mm}

$$\bigg(\int_{\sqrt{t}}^t\Big|\sum_{\substack{\ell\leq u\\P^-(\ell)>y}}\Lambda(\ell)\chi(\ell)(\log \ell)^k\Big|\frac{\dee u}{u}\bigg)^2\leq\frac{\log t}{2}\int_{\sqrt{t}}^t\Big|\sum_{\substack{\ell\leq u\\P^-(\ell)>y}}\Lambda(\ell)\chi(\ell)(\log \ell)^k\Big|^2\frac{\dee u}{u},$$

\vspace{1mm}

\noindent
by the Cauchy-Schwarz inequality, we conclude that

\vspace{-2mm}

\begin{align*}
&\bigg(\sum_{\substack{\chi\mod{q}\\\chi\neq \chi_0, \psi}}\int_{\frac{\sqrt{x}}{2}}^x\Big|\sum_{\substack{\ell\leq t\\P^-(\ell)>y}}\Lambda(\ell)\chi(\ell)\Big|^2\frac{\dee t}{t^3}\bigg)^{1/2}\\
&\ll\frac{\phi(q)^{1/2}}{x^{1/4}}+M^k(\log x)^{-k}\bigg(\sum_{\substack{\chi\mod{q}\\\chi\neq \chi_0, \psi}}\int_{\frac{\sqrt{x}}{2}}^x\Big|\sum_{\substack{\ell\leq t\\P^-(\ell)>y}}\Lambda(\ell)\chi(\ell)(\log \ell)^k\Big|^2\frac{\dee t}{t^3}\bigg)^{1/2}\nonumber\\
&+M^k(\log x)^{-k}\bigg(\sum_{\substack{\chi\mod{q}\\\chi\neq \chi_0, \psi}}\int_{\frac{\sqrt{x}}{2}}^x\bigg(\int_{\sqrt{t}}^t\Big|\sum_{\substack{\ell\leq u\\P^-(\ell)>y}}\Lambda(\ell)\chi(\ell)(\log \ell)^k\Big|^2\frac{\dee u}{u}\bigg)\frac{\dee t}{t^3}\bigg)^{1/2}.\nonumber
\end{align*}

\vspace{1mm} 

\noindent
We use Fubini's theorem to interchange the order of integration, and so the double integral above equals 

\vspace{-2mm}

\begin{align*}
\int_{\frac{\sqrt[4]{x}}{\sqrt{2}}}^x\Big|\sum_{\substack{\ell\leq u\\P^-(\ell)>y}}\Lambda(\ell)\chi(\ell)(\log \ell)^k\Big|^2\bigg(\int_u^{u^2}\frac{\dee t}{t^3}\bigg)\frac{\dee u}{u}\leq\int_{\frac{\sqrt[4]{x}}{2}}^x\Big|\sum_{\substack{\ell\leq u\\P^-(\ell)>y}}\Lambda(\ell)\chi(\ell)(\log \ell)^k\Big|^2\frac{\dee u}{u^3}.
\end{align*}

\noindent
Hence,

\vspace{-2mm}

\begin{align*}
&\bigg(\sum_{\substack{\chi\mod{q}\\\chi\neq \chi_0, \psi}}\int_{\frac{\sqrt{x}}{2}}^x\Big|\sum_{\substack{\ell\leq t\\P^-(\ell)>y}}\Lambda(\ell)\chi(\ell)\Big|^2\frac{\dee t}{t^3}\bigg)^{1/2}\\
&\ll \frac{\phi(q)^{1/2}}{x^{1/4}} +M^k(\log x)^{-k}\bigg(\sum_{\substack{\chi\mod{q}\\\chi\neq \chi_0, \psi}}\int_{\frac{\sqrt[4]{x}}{2}}^x\Big|\sum_{\substack{\ell\leq t\\P^-(\ell)>y}}\Lambda(\ell)\chi(\ell)(\log \ell)^k\Big|^2\frac{\dee t}{t^3}\bigg)^{1/2}\nonumber\\
&\ll \frac{\phi(q)^{1/2}}{x^{1/4}} +M^k(\log x)^{-k}\bigg(\sum_{\substack{\chi\mod{q}\\\chi\neq \chi_0, \psi}}\int_{\frac{\sqrt[4]{x}}{2}}^x\Big|\sum_{\substack{\ell\leq t\\P^-(\ell)>y}}\Lambda(\ell)\chi(\ell)(\log \ell)^k\Big|^2\frac{\dee t}{t^{3+2/\log x}}\bigg)^{1/2}.\nonumber
\end{align*}

\noindent
Since Parseval's theorem for Dirichlet series guarantees that

$$\int_1^{\infty}\Big|\sum_{\substack{\ell\leq u\\P^-(\ell)>y}}\Lambda(\ell)\chi(\ell)(\log \ell)^k\Big|^2\frac{\dee u}{u^{3+2/\log x}}=\frac{1}{2\pi}\int_{\R}\Big|\bigg(\frac{L_y^{'}}{L_y}\bigg)^{(k)}(c+it,\chi)\Big|^2\frac{\dee t}{c^2+t^2},$$

\noindent
with $c:=1+1/\log x$, we arrive at the following bound:

\begin{align}\label{afP}
&\bigg(\sum_{\substack{\chi\mod{q}\\\chi\neq \chi_0, \psi}}\int_{\frac{\sqrt{x}}{2}}^x\Big|\sum_{\substack{\ell\leq t\\P^-(\ell)>y}}\Lambda(\ell)\chi(\ell)\Big|^2\frac{\dee t}{t^3}\bigg)^{1/2}\\
&\ll \frac{\phi(q)^{1/2}}{x^{1/4}} +M^k(\log x)^{-k}\bigg(\sum_{\substack{\chi\mod{q}\\\chi\neq\chi_0, \psi}}\int_{\R}\Big|\bigg(\frac{L_y^{'}}{L_y}\bigg)^{(k)}(c+it,\chi)\Big|^2\frac{\dee t}{c^2+t^2}\bigg)^{1/2}.\nonumber
\end{align}

We now restrict our attention to the integrals involving the derivatives of $L_y^{'}/L_y$. We are going to split these integrals into two parts. In the first parts, we will be integrating over $|t|\leq T$. We will bound these parts by mainly using the results of Section 3. For the remaining parts, where we integrate over the range $|t|>T$, we will use Lemma \ref{mvt}. We start with the integrals over $|t|>T$ first. 

It is true that

\begin{align*}
\sum_{\chi\mod{q}}\bigg|\bigg(\frac{L_y^{'}}{L_y}\bigg)^{(k)}(c+it,\chi)\bigg|^2=\phi(q)\sum_{b\in(\Z/q\Z)^*}\bigg|\sum_{\substack{n\equiv b\mod{q}\\P^-(n)>y}}\frac{\Lambda(n)(\log n)^k}{n^{c+it}}\bigg|^2.
\end{align*}

\noindent
This may be proven by opening the square of the left hand-side and by using the orthogonality relations of the characters modulo $q$. So, now one can use Lemma \ref{mvt} to infer that

\begin{align}\label{smt}
&\sum_{\substack{\chi\mod{q}\\\chi\neq\chi_0, \psi}}\int_{|t|>T}\Big|\bigg(\frac{L_y^{'}}{L_y}\bigg)^{(k)}(c+it,\chi)\Big|^2\frac{\dee t}{c^2+t^2}\\
&\leq\phi(q)\sum_{b\in(\Z/q\Z)^*}\int_{|t|>T}\bigg|\sum_{\substack{n\equiv b\mod{q}\\P^-(n)>y}}\frac{\Lambda(n)(\log n)^k}{n^{c+it}}\bigg|^2\frac{\dee t}{t^2}\ll \frac{(8k)^{2k}(\log x)^{2k+1}}{T},\nonumber\end{align}

Our treatment for the integrals corresponding to the large values of $t$ is complete and we turn our focus to the integrals whose range of integration is $|t|\leq T$. For a character $\chi\notin\{\chi_0,\psi\}$, using Lemma \ref{lasl}, it follows that

\begin{align}\label{afld}
&\int_{|t|\leq T}\Big|\bigg(\frac{L_y^{'}}{L_y}\bigg)^{(k)}(c+it,\chi)\Big|^2\frac{\dee t}{c^2+t^2}\\
&\leq4^{k+1}((k+1)!)^2\sum_{j=1}^{k+1}(j!)^{-\frac{2(k+1)}{j}}\int_{|t|\leq T}\Big|\frac{L_y^{(j)}}{L_y}(c+it,\chi)\Big|^{\frac{2(k+1)}{j}}\frac{\dee t}{c^2+t^2}.\nonumber
\end{align}

\vspace{3mm}

\noindent
Since $\chi\notin\{\chi_0,\psi\}$ and $y=(10q)^{100}V_T>qV_t$ when $|t|\leq T$, a proper combination of Theorems \ref{lb} and \ref{ldnS} implies that $\abs{L_y^{-1}(c+it,\chi)}\ll 1$. Moreover, $|L_y^{(j)}(c+it,\chi)|\ll j!(C\log y)^j$ for all $j\in\{1,\ldots,k+1\}$, as can be seen from Theorem \ref{ub}. Hence, (\ref{afld}) gives

\begin{align}\label{1rel}
&\sum_{\substack{\chi\mod{q}\\\chi\neq\chi_0,\psi}}\int_{|t|\leq T}\Big|\bigg(\frac{L_y^{'}}{L_y}\bigg)^{(k)}(c+it,\chi)\Big|^2\frac{\dee t}{c^2+t^2}\\
&\ll (2C)^{2k}((k+1)!)^2\sum_{j=1}^{k+1}\frac{(\log y)^{2(k+1-j)}}{(j!)^2}\sum_{\substack{\chi\mod{q}\\\chi\neq\chi_0,\psi}}\int_{\R}\abs{L_y^{(j)}(c+it,\chi)}^2\frac{\dee t}{c^2+t^2}.\nonumber
\end{align}

\vspace{2mm}

\noindent
From an application of Parseval's theorem, it follows that

\begin{align}\label{lP}
\sum_{\substack{\chi\mod{q}\\\chi\neq\chi_0, \psi}}\int_{\R}\abs{L_y^{(j)}(c+it,\chi)}^2\frac{\dee t}{c^2+t^2}\leq\int_y^{\infty}\sum_{\substack{\chi\mod{q}\\\chi\neq\chi_0}}\Big|\sum_{\substack{n\leq u\\P^-(n)>y}}\chi(n)(\log n)^j\Big|^2\frac{\dee u}{u^3}.
\end{align}

\noindent
For $b\in(\Z/q\Z)^*$ and $u\geq y$, Lemmas \ref{chlem} and \ref{logap} yield that

\begin{align}\label{1''}
\sum_{\substack{n\leq u\\n\equiv b \mod{q}\\P^-(n)>y}}(\log n)^j-\frac{1}{\phi(q)}\sum_{\substack{n\leq u\\P^-(n)>y}}\chi_0(n)(\log n)^j\ll\frac{(\log u)^ju^{1-c_1/\log y}}{\phi(q)\log y}.
\end{align}

\noindent
Since $\int_y^v(\log t)^j\dee t\leq v(\log v)^j$ for $v\geq y$, these lemmas also imply that

\begin{align}\label{2''}
\sum_{\substack{n\leq u\\n\equiv b \mod{q}\\P^-(n)>y}}(\log n)^j+\frac{1}{\phi(q)}\sum_{\substack{n\leq u\\P^-(n)>y}}\chi_0(n)(\log n)^j\ll\frac{u(\log u)^j}{\phi(q)\log y}.
\end{align}

\noindent
We combine (\ref{1''}) and (\ref{2''}) with the elementary identity $w^2-z^2=(w-z)(w+z)$ and infer that 

\begin{align*}
&\sum_{\substack{\chi\mod{q}\\\chi\neq \chi_0}}\Big|\sum_{\substack{n\leq u\\P^-(n)>y}}\chi(n)(\log n)^j\Big|^2\\
&=\phi(q)\sum_{b\in(\Z/q\Z)^*}\bigg\{\Big(\sum_{\substack{n\leq u\\n\equiv b \mod{q}\\P^-(n)>y}}(\log n)^j\Big)^2-\Big(\frac{1}{\phi(q)}\sum_{\substack{n\leq u\\P^-(n)>y}}\chi_0(n)(\log n)^j\Big)^2\bigg\}\\
&\ll\frac{(\log u)^{2j}u^{2-c_1/\log y}}{(\log y)^2},
\end{align*}

\vspace{3mm}

\noindent
for all $j\in\{1,\ldots,k+1\}$. We plug this estimate into (\ref{lP}) and obtain

\vspace{1mm}

\begin{align*}
\sum_{\substack{\chi\mod{q}\\\chi\neq\chi_0, \psi}}\int_{\R}\abs{L_y^{(j)}(c+it,\chi)}^2\frac{\dee t}{c^2+t^2}&\ll\frac{1}{(\log y)^2}\int_y^{\infty}(\log u)^{2j}u^{-1-c_1/\log y}\dee u\\
&=c_1^{-2j-1}(\log y)^{2j-1}\Gamma(2j+1)\leq c_1^{-3}(2j)!(\log y)^{2j-1},
\end{align*}

\vspace{2mm}

\noindent
for $j\in\{1,\ldots,k+1\}$ and $\Gamma$ being the Gamma function. With this last bound, (\ref{1rel}) turns into

\begin{align*}
\sum_{\substack{\chi\mod{q}\\\chi\neq\chi_0,\psi}}\int_{|t|\leq T}\Big|\bigg(\frac{L_y^{'}}{L_y}\bigg)^{(k)}(c+it,\chi)\Big|^2\frac{\dee t}{c^2+t^2}\ll (2C)^{2k}((k+1)!)^2(\log y)^{2k+1}\sum_{j=1}^{k+1}\binom{2j}{j}.
\end{align*}

\noindent
But, $\binom{2j}{j}\leq \sum_{\ell=0}^j\binom{2j}{\ell}=4^{j}$, and so

$$\sum_{j=1}^{k+1}\binom{2j}{j}\leq 4^{k+1}\sum_{j=0}^{k}4^{-j}<4^{k+1}\sum_{j\geq 0}4^{-j}\ll 4^k.$$

\vspace{2mm}

\noindent
Therefore, we complete the estimation of the integrals over $|t|\leq T$ by arriving at the bound

\begin{align}\label{lart}
\sum_{\substack{\chi\mod{q}\\\chi\neq\chi_0,\psi}}\int_{|t|\leq T}\Big|\bigg(\frac{L_y^{'}}{L_y}\bigg)^{(k)}(c+it,\chi)\Big|^2\frac{\dee t}{c^2+t^2}\ll (4C)^{2k}((k+1)!)^2(\log y)^{2k+1}.
\end{align}

Putting (\ref{smt}) and (\ref{lart}) together, we deduce that

\begin{align*}
\sum_{\substack{\chi\mod{q}\\\chi\neq\chi_0,\psi}}\int_{\R}\Big|\bigg(\frac{L_y^{'}}{L_y}\bigg)^{(k)}(c+it,\chi)\Big|^2\frac{\dee t}{c^2+t^2}\ll (4C)^{2k}((k+1)!)^2(\log y)^{2k+1}+\frac{(8k)^{2k}(\log x)^{2k+1}}{T}.
\end{align*}

\vspace{2mm}

\noindent
We insert this into (\ref{afP}) and conclude that

\begin{align}\label{afc}
&\bigg(\sum_{\substack{\chi\mod{q}\\\chi\neq \chi_0, \psi}}\int_{\frac{\sqrt{x}}{2}}^x\Big|\sum_{\substack{\ell\leq t\\P^-(\ell)>y}}\Lambda(\ell)\chi(\ell)\Big|^2\frac{\dee t}{t^3}\bigg)^{1/2}\\
&\ll \frac{\phi(q)^{1/2}}{x^{1/4}} +(4MC)^k(k+1)!\frac{(\log y)^{k+\frac{1}{2}}}{(\log x)^k}+\frac{(8Mk)^k\sqrt{\log x}}{\sqrt{T}}.\nonumber
\end{align}

\vspace{3mm}

\noindent
In turn, (\ref{afc}), when used in (\ref{bP}), implies that

\begin{align}\label{alend}
\frac{1}{D}\sum_{\substack{1<m\leq\sqrt{x}\\P^-(m)>y}}\int_{x-D}^x\Delta\Big(\frac{t}{m},y;q,a\overline{m}\Big)\dee t\ll_{\delta}\frac{x^{1+2\delta/\log y}}{\phi(q)}\bigg\{\bigg(\frac{M'\ell\log y}{\log x}\bigg)^{\!\!\ell}+\frac{(8M\ell)^{\ell}}{\sqrt{T}}\bigg\},
\end{align}

\vspace{2mm}

\noindent
where $\ell=k+1$ and $M'=4MC$. Note that we omitted the term stemming from $\phi(q)^{1/2}x^{-1/4}$, since we are working with a sufficiently small $\delta>0$ and a $x\geq q^B$ for some sufficiently large $B$. Now, we combine (\ref{recres}) with (\ref{alend}) and obtain that

\begin{align}\label{forend}
\Delta(x,y;q,a)\ll_{\delta} \frac{x^{1+2\delta/\log y}}{\phi(q)}\bigg\{\bigg(\frac{M'\ell\log y}{\log x}\bigg)^{\!\!\ell}+\frac{(8M\ell)^{\ell}}{\sqrt{T}}\bigg\}+\frac{x^{1-\delta/(2\log y)}}{\phi(q)}.
\end{align}

\vspace{2mm} 

Since $M$ can be sufficiently large, we can now choose $\delta$ to be sufficiently small as $\delta=1/(5eM')$. Then, for 

$$\ell=\flbgg{\frac{\log x}{eM'\log y}},$$ 

\noindent
it follows that

$$\bigg(\frac{M'\ell\log y}{\log x}\bigg)^{\!\!\ell}\leq \frac{\log x}{M'\log y}x^{-1/(eM'\log y)}\ll x^{-1/(2eM'\log y)},$$

\vspace{3mm}

\noindent
and so, with the above choice of $\delta$, we deduce that

\begin{align}\label{opt1}
\frac{x^{1+2\delta/(\log y)}}{\phi(q)}\bigg(\frac{M'\ell\log y}{\log x}\bigg)^{\!\!\ell}\ll \frac{x^{1-\delta/(2\log y)}}{\phi(q)}.
\end{align}

\vspace{3mm}

\noindent 
Now, since $T=\exp\{2L(\log x)^{3/5}(\log\log x)^{2/5}\}$, we observe that $\log V_T\asymp (\log x)^{2/5}(\log \log x)^{3/5}$, and this implies that $\ell<\log x/\log y<\log x/\log V_T\ll (\log x)^{3/5}(\log \log x)^{-3/5}$. Moreover, $\log \ell<\log\log x$. Consequently, there exist some positive constants $c_3$ and $c_4$ such that
\begin{align}\label{opt2}
x^{\frac{2\delta}{\log y}}\frac{(8M\ell)^{\ell}}{\sqrt{T}}&= \exp\bigg\{\frac{2\delta\log x}{\log y}+\ell\log \ell+\log(8M)\ell-L(\log x)^{3/5}(\log\log x)^{2/5}\bigg\}\nonumber\\
&\leq\exp\{c_3(\log x)^{3/5}(\log\log x)^{1-3/5}-L(\log x)^{3/5}(\log\log x)^{2/5}\}\nonumber\\
&\leq\exp\{-c_4(\log x)^{3/5}(\log\log x)^{2/5}\},
\end{align}

\vspace{2mm}

\noindent
because the constant $L$ in the definition of $T$ is sufficiently large. 

With the selection of $\ell$ that we made, we have the estimates (\ref{opt1}) and (\ref{opt2}) and then (\ref{forend}) becomes
\begin{align}\label{toendit}
\Delta(x,y;q,a)&\ll\frac{x^{1-\delta/(2\log y)}}{\phi(q)}+\frac{xe^{-c_4(\log x)^{3/5}(\log\log x)^{2/5}}}{\phi(q)}\\
&\ll\frac{x^{1-c_5/(\log(2q))}}{\phi(q)}+\frac{xe^{-c_5(\log x)^{3/5}(\log\log x)^{-3/5}}}{\phi(q)},
\end{align}
\noindent
for some $c_5>0$. The first term of the second line is for the range where $\log(10q)\geq \log V_T$, whereas the second term covers the range $\log(10q)\leq \log V_T$. The proof of the theorem is almost complete. It only remains to observe that $\Delta(x,y;q,a)$ does not differ much from

$$\sum_{\substack{n\leq x\\n\equiv a\mod{q}}}\Lambda(n)-\frac{x}{\phi(q)}-\frac{\psi(a)}{\phi(q)}\sum_{n\leq x}\Lambda(n)\psi(n).$$

\vspace{2mm}

First,

$$\sum_{\substack{n\leq x\\n\equiv a\mod{q} \\P^-(n)>y}}\Lambda(n)=\sum_{\substack{n\leq x\\n\equiv a\mod{q}}}\Lambda(n)+O\bigg(\frac{y\log x}{\log y}\bigg),$$

\noindent
because

\begin{align*}
0\leq \sum_{\substack{n\leq x\\n\equiv a\mod{q}}}\Lambda(n)-\sum_{\substack{n\leq x\\n\equiv a\mod{q} \\P^-(n)>y}}\Lambda(n)\leq\sum_{\substack{p^{\nu}\leq x\\p\leq y}}\Lambda(p^{\nu})\leq \pi(y)\log x\ll\frac{y\log x}{\log y}.
\end{align*}

\noindent
Similarly,

$$\sum_{\substack{n\leq x\\P^-(n)>y}}\Lambda(n)\chi(n)=\sum_{n\leq x}\Lambda(n)\chi(n)+O\bigg(\frac{y\log x}{\log y}\bigg)$$

\vspace{1mm}

\noindent
for $\chi\in\{\chi_0,\psi\}$. In addition, there exists a constant $c_6>0$ such that (this can also follow from the pretentious methods pursued in \cite{ppnt})

$$\sum_{n\leq x}\Lambda(n)\chi_0(n)=x+O(x\exp\{-c_6(\log x)^{3/5}(\log\log x)^{-1/5}\}),$$

\vspace{1mm}

\noindent
and so, when $x\geq q^B$ for a large $B>0$, with the $y$ that we have chosen, we conclude that

\begin{align*}
\sum_{\substack{n\leq x\\n\equiv a\mod{q}}}\Lambda(n)-\frac{x}{\phi(q)}-\frac{\psi(a)}{\phi(q)}\sum_{n\leq x}\Lambda(n)\psi(n)-\Delta(x,y;q,a)\ll\frac{xe^{-c_6(\log x)^{3/5}(\log \log x)^{-1/5}}}{\phi(q)}
\end{align*}

\vspace{3mm}

\noindent
In virtue of (\ref{toendit}), the theorem has now been proven. 


\bibliographystyle{alpha}

\end{document}